\documentclass[11pt]{article}
\usepackage{amssymb,amsfonts}
\usepackage{CJK}
\def\<{\langle}
\def\>{\rangle}
\def\a{\alpha}

\def\c{\cdot}
\def\ci{\circ}

\def\D{\Delta}
\def\e{\eta}

\def\m{\mu}

\def\lr{\longrightarrow}

\def\o{\otimes}

\def\lm{\longmapsto}
\def\ol{\overline}

\def\ul{\underline}
\def\v{\varepsilon}

\def\R{\mathcal{R}}
\def\wt{\widetilde}
\date{}
\begin{document}
\renewcommand{\baselinestretch}{1.2}
\renewcommand{\arraystretch}{1.0}
\title{\bf Transmutation Theory and Quantization Approach for Quantum Groupoids}
\date{}
\author {{\bf Xuan Zhou\footnote {Corresponding author. Mathematics and Information Technology School,
Jiangsu Second Normal University, Nanjing, Jiangsu 210013, China. E-mail: zhouxuanseu@126.com}\ ,
Tao Yang\footnote {College of Science, Nanjing Agricultural University, Nanjing, Jiangsu 210095, China.}\ ,
 }\\
}
\maketitle
\begin{center}
\begin{minipage}{12.cm}

 {\bf Abstract:} Let $H$ and $L$ be quantum groupoids.
 If $H$ has a quasitriangular structure, then we show that $L$ induces a Hopf algebra $C_{L}(L_s)$ in the category $_{H}\mathcal{M}$,
 which generalizes the transmutation theory introduced by Majid.
 Furthermore, if $H$ is commutative, we can construct a Hopf algebra $C_H(H_s)_F$ in the category $_H\mathcal{M}_F$
 for a weak invertible unit 2-cocycle $F$,
 which generalizes the results in \cite{D83}.
 Finally, we consider the relation between two Hopf algebras: $C_H(H_s)_F$ and $C_{\widetilde H}(\widetilde{H}_s)$, and obtain that
 they are isomorphic as objects in the category $_{\widetilde H}\mathcal{M}$, where $(\widetilde H, \widetilde{\mathcal{R}})$ is a new quasitriangular
 quantum groupoid induced by $(H, \mathcal{R})$.

\vskip 0.5cm


\end{minipage}
\end{center}

\section*{0. Introduction}
\def\theequation{1. \arabic{equation}}
\setcounter{equation} {0} \hskip\parindent

 Hopf algebras in braided monoidal categories (i.e., braided groups), introduced by Majid in the study
 of Tannaka-Krein-type reconstruction theory \cite{M93,M91}, play an important role in mathematics and physics.
 The method (so-called transmutation) that every quantum group gives rise to a Hopf algebra in a braided category
 was provided in \cite{M93}. And the resulting braided group
 is braided-cocommutative in the representation category of the quantum group.
 \\

 And Hopf algebras in symmetric monoidal categories came
 from the deformation-quantization of triangular solutions of the classical Yang-Baxter
 equations \cite{GRZ92}. They are called $S$-Hopf algebras, which are obtained
 by means of an element $F$ constructed by Drinfel'd in \cite{D83}. Moreover, Gurevich and Majid
 \cite{GM94} pointed out this method can be used quite generally and defined the $S$-Hopf algebra
 $H_F$ for any pair $(H, F)$, where $H$ is a cocommutative Hopf algebra and $F$ satisfies a cocycle
 condition.
\\

 As one of main generalizations of Hopf algebras, quantum groupoids (i.e., weak Hopf algebras) were
 defined by B\"ohm et. al \cite{BNS99}. The motivation for studying quantum groupoids comes from quantum
 field theories and operator algebras. And the representation theory of quantum groupoids provides examples
 of monoidal categories that can be used for constructing invariants of links and $3$-manifolds \cite{NTV03}.
 Although many results of classical Hopf algebra theory could be generalized to quantum groupoids, the
 processes of proof and derivation are very complex.
\\

 So, there is a natural question: whether braided groups can be constructed from quantum groupoids, that is,
 whether the transmutation theory and quantization approach can be
 generalized to the quantum groupoids setting.
 Let $(H, \R$) be a quasitriangular quantum groupiod, $L$ a quantum groupoid, and there is a quantum groupoid
 map between them. Although we could not obtain $L$
 is a Hopf algebra in the category $_H\mathcal{M}$ of left
 $H$-modules, we found that the centralizer algebra $C_{L}(L_s)$ meets our needs.
 Using the similar technique to consider the quantization approach, we obtain a Hopf algebra $C_H(H_s)_F$
 in the category $_H{\mathcal{M}}_F$, where $H$ satisfies the cocommutative condition and $F$ is
 a weak invertible unit 2-cocycle.
 \\

Furthermore, the relationship between the two methods of constructing braided groups over quantum groupoids
is another purpose of this paper. The twisting construction, introduced by Drinfel'd \cite{D90}, provides a way to
obtain new quasitriangular Hopf algebras from the given ones and the special elements.
In 2008, Chen \cite{C08} generalized the twisting theory to the quantum groupiods setting. Applying
our theory to the new quasitriangular quantum groupoid $(\wt H, \wt\R)$, we get
a Hopf algebra $C_{\wt H}(\wt{H}_s)$ in the category $_{\wt H}\mathcal{M}$ of the left
$\wt H$-modules and obtain the isomorphism theorem between $C_{\wt H}(\wt{H}_s)$ and $C_H(H_s)_F$.
\\

The paper  is  organized  as  follows.
\\

  In Section 1, we recall some basic notions and results
  for quantum groupoids and braided monoidal categories. In
  Section 2, we obtain a Hopf algebra $C_L(L_s)$ in the category $_H\mathcal{M}$ of
  left $H$-modules through the quasitriangular structure of $H$. In Section 3,
  we construct a Hopf algebra $C_H(H_s)_F$ in the category $_H\mathcal{M}_F$, for a weak invertible unit 2-cocycle $F$.
  In Section 4, by the weak twisting theory, we discuss the relationship between the above two methods.

\section*{1. Preliminaries}
\def\theequation{1. \arabic{equation}}
\setcounter{equation} {0} \hskip\parindent

Throughout the paper, we let $k$ be a fixed
 field, and use the Sweedler formal sum notation for the comultiplication
 $\D$ over a coalgebra $C$ \cite{Sw69}, that is, $\Delta(c)=c_{1}\o c_{2}$ for any $c\in C$.

\vskip 0.5cm
 {\bf 1.1 Basic definitions and properties about quantum groupoids}
\vskip 0.5cm

 Recall from B$\ddot{\mbox o}$hm et al.\cite{BNS99} that
 a weak bialgebra $H$
 is an algebra $(H, \m, \e)$ and a coalgebra $(H, \D, \v)$
 such that $\Delta(hl)=\Delta(h)\Delta(l)$ and
\begin{eqnarray*}
&&\Delta^{2}(1)=1_{1}\o 1_{2}1'_{1}\o 1'_{2} =1_{1}\o
1'_{1}1_{2}\o 1'_{2}=1_{1}\o 1_{2}\o 1_{3},\\
&&\v(hgl)=\v(hg_{1})\v(g_{2}l)=\v(hg_{2})\v(g_{1}l),
\end{eqnarray*}
 for all $h, g, l\in H$.
\\

Moreover, a weak bialgebra $H$ is called a quantum groupoid if
there exists a $k$-linear map $S: H\lr H$, satisfying the following
conditions:
$$
S* id =\varepsilon_{s}, \ \  id *S=\varepsilon_{t}, \ \ S*id*S=S,
$$
where
 $*$ is the usual convolution product, and the idempotent maps $\v_t, \v_s$
 are defined by
$$
\v_{t}(h)=\v(1_{1}h)1_{2}, \ \ \v_{s}(h)=1_{1}\v(h1_{2}), \ \ \mbox{for all} \ \ h\in H.
$$

In this case, $S$ is called the antipode, $\v_t$ is called the target map, and $\v_s$ is called
the source map, respectively. If antipode $S$ exists, then it is unique and  $S$ is an
 anti-algebra and anti-coalgebra morphism. We will always assume that $S$ is
 bijective. If $H$ is a finite dimensional quantum groupoid,
 then $S$ is automatically bijective.
\\

The target space $H_{t}$ and source space $H_{s}$ are the images of $\v_t$ and $\v_s$, respectively,
which can be described as follows:
\begin{eqnarray*}
&& H_t=\{h\in H\mid \v _t(h)=h\}=\{h\in H\mid \Delta
(h)=1_{1}h\o 1_{2}=h1_{1}\o 1_{2}\},\\
&& H_s=\{h\in H\mid \v _s(h)=h\}=\{h\in H\mid \Delta (h)=1_{1}\o
 h1_{2}=1_{1}\o 1_{2}h\}.
\end{eqnarray*}
It is easy to see that $H_{t}$ and $H_{s}$ are subalgebras of $H$.

Similarly, we have
 \begin{eqnarray*}
 \ol{\varepsilon}_{s}(h)=\varepsilon(h1_{1})1_{2};\quad
\ol{\varepsilon}_{t}(h)=1_{1}\varepsilon(1_{2}h).
 \end{eqnarray*}

 Let $H$ be a quantum groupoid. Recall from Chen \cite{C08}, we say an element $F\in \D(1)(H\o H)\D^{cop}(1)$
is a weak invertible unit 2-cocycle if it has the following properties:

(1) There exists an element $F^{-1}\in  \D^{cop}(1)(H\o H)\D(1)$, such that
$$
FF^{-1}=\D(1), \ \ \ F^{-1}F=\D^{cop}(1).
$$

(2)$((\D\o id)F)F_{12}=((id\o\D)F)F_{23}$, where $F_{23}=1\o F$, $F_{12}=F\o 1$.

(3) For all $y\in H_s$, and $z\in H_t$, the following equations hold:
\begin{tabbing}
lllllllllllll \= 22222222222222222222222222222222 \= \kill
\> $(1\o y)F=F(y\o 1)$ , \> $(z\o1)F=F(1\o z)$ ; \\
\> $F^{-1}(1\o y)=(y\o 1)F^{-1}$ , \> $F^{-1}(z\o 1)=(1\o z)F^{-1}$ ; \\
\> $(1\o y)F^{-1}=(S^{-1}(y)\o 1)F^{-1}$ , \> $F(z\o 1)=F(1\o S^{-1}(z))$ .
\end{tabbing}

We list some equivalent forms of $((\D\o id)F)F_{12}=((id\o\D)F)F_{23}$ for convenience
as follows:
$$
\left\{%
\begin{array}{l}
F'^{-(1)}\o F'^{(1)}F'^{-(2)}\o F'^{(2)}=F^{-(1)}F^{(1)}_1\o F^{-(2)}_1F^{(1)}_2\o F^{-(2)}_2F^{(2)}, \\
F'^{(1)}\o F'^{(2)}F'^{-(1)}\o F'^{-(2)}=F^{-(1)}_1F^{(1)}\o F^{-(1)}_2F^{(2)}_1\o F^{-(2)}F^{(2)}_2, \\
F^{-(1)}\o F'^{-(1)}F^{-(2)}_1\o F'^{-(2)}F^{-(2)}_2=F'^{-(1)}F^{-(1)}_1\o F'^{-(2)}F^{-(1)}_2\o F^{-(2)},
\end{array}%
\right.
$$
where $F'=F$ and $F'^{-1}=F^{-1}$.
\\

We recall from Nikshych et al. \cite{NTV03} that
a quasitriangular quantum groupoid is a pair $(H,\R)$, where $H$
is a quantum groupoid, and $\R=R^{(1)}\otimes
R^{(2)}\in \Delta^{cop}(1)(H\otimes_{k}H)\Delta(1)$, satisfying the
following conditions:
$$
(id\otimes \Delta)\R=\R_{13}\R_{12}, \ \ (\Delta\otimes
id)\R=\R_{13}\R_{23},
$$
$$
\Delta^{cop}(h)\R=\R\Delta(h), \ \ \ \mbox{for all}\ \ h\in H,
$$
 where $\R_{12}=\R\otimes
1$, $\R_{23}=1\otimes \R$ etc., as usual, and such that there exists
$\R^{-1}=R^{-(1)}\otimes
R^{-(2)}\in \Delta(1)(H\otimes_{k}H)\Delta^{cop}(1)$ with

$$
\R\R^{-1}=\Delta^{cop}(1),  \ \ \ \R^{-1}\R=\Delta(1).
$$

Next, we summarize some important properties of quasitriangular quantum groupoids from
\cite{NTV03} as follows, which is needed for our calculation through the paper.
\\

Let $(H, \R)$ be a quasitriangular quantum groupoid. Then
 we have the following equations:
\begin{tabbing}
lllllllll \= 222222222222222222222222222222222222 \= \kill
\> $(1\o z)\R=\R(z\o 1)$ , \> $(y\o1)\R=\R(1\o y)$ , \\
\> $(z\o1)\R=(1\o s(z))\R$ , \> $(1\o y)\R=(s(y)\o 1)\R$ , \\
\> $\R(1\o z)=\R(S^{-1}(z)\o1)$ , \> $\R(y\o1)=\R(1\o S^{-1}(y))$ ,\\
\> $(\v_s\o id)(\R)=\D(1)$ , \> $(id\o\v_s)(\R)=(S\o
id)\D^{cop}(1)$ , \\
\> $(\v_t\o id)(\R)=\D^{cop}(1)$ , \> $(id\o\v_t)(\R)=(S\o
id)\D(1)$ , \\
\> $(S\o id)(\R)=(id\o S^{-1})(\R)=\R^{-1}$ , \> $(S\o
S)(\R)=\R$ ,
\end{tabbing}
for all $y\in H_s, z\in H_t$.
\\

For all $h\in H$, we have $S^2(h)=uhu^{-1}$, where
$u=S(R^{(2)})R^{(1)}$ is an invertible element of $H$ such that
$$
u^{-1}=R^{(2)}S^2(R^{(1)}),\ \ \ \ \D(u)=\R^{-1}\R^{-1}_{21}(u\o
u).
$$

\vskip 0.5cm
 {\bf 1.2  Some definitions in braided monoidal categories}
\vskip 0.5cm

We denote a braided monoidal category $\mathcal{C}$ as $({\mathcal C}, \o , I, a, l, r, \Psi)$,
where $\mathcal{C}$ is category with tensor product $\o$, unit $I$,
associator $a$, left and right unit constraints $l$ and $r$.
The braiding is denoted by $\Psi$. Then we review the definition of Hopf algebras in category $\mathcal{C}$,
see \cite{ARFLV03}, \cite{Ka95}, \cite{M02}.
\\

An algebra in ${\mathcal C}$ is a triple $A=(A, \e,\m)$ where
$A$ is an object in ${\mathcal C}$ and $\e: I\lr A$ and
$\m: A \o A \lr A$ are morphisms in
${\mathcal C}$ such that $\m\circ (id\o \e)=id=\m
\circ (\e\o id)$, and $ \m \circ (id\o \m)=\m \circ
(\m\o id)$.
\\

 A coalgebra in ${\mathcal C}$ is a triple $C=(C, \v, \D)$ where
$C$ is an object in ${\mathcal C}$ and $\v: C \lr I$ and
$\D: C\lr C \o C$ are morphisms in
${\mathcal C}$ such that $(id\o \v)\circ \D=id=(\v\o
id)\circ \D$, and $(id\o \D)\circ \D=(\D\o id)\circ \D$.
\\

A bialgebra $H$ is an object in $\mathcal{C}$ with an algebra structure $(H, \e,\m)$
and a coalgebra structure $(H, \v, \D)$ such that the following axioms hold:
\begin{eqnarray*}
&&\D\circ \m=(\m\o \m)\circ (id\o \Psi\o
id)\circ(\D\o \D),\\
&&\v\circ \m= \v \o \v, \ \ \ \ \  \D(1)=1\o 1.
\end{eqnarray*}
A bialgebra $(H, \m, \e, \D, \v)$ in ${\mathcal C}$ is
called a Hopf algebra in ${\mathcal C}$ if there exists a
morphism $S: H\lr H$ in ${\mathcal C}$ satisfying
\begin{eqnarray*}
\m\circ (S\o id)\circ\D=\v\circ \e=\m\circ (id\o
S)\circ \D.
\end{eqnarray*}

\section*{2. Transmutation theory for quantum groupoids}
\def\theequation{1. \arabic{equation}}
\setcounter{equation} {0} \hskip\parindent

 Let $H$ be a quantum groupoid with an antipode $S$, $_{H}\mathcal{M}$
denote the category of left $H$-modules. In this section, we suppose that $H$ endowed with a quasitriangular structure $\R$.
Then we get a monoidal category $(_{H}\mathcal{M}, \o_t, H_t, a, l, r)$ \cite{NTV03}, where $l$ and $r$ are defined by
$$
l(z\o_t v)=z\c v, \ \ r(v\o_t z)=S^{-1}(z)\c v, \ \ \mbox{for any}\ \  z\in H_t, v\in V\in {_H\mathcal{M}}.
$$
 Furthermore, the category $(_{H}\mathcal{M} ,
 \o_{t}, H_{t}, a, l, r)$ is braided,
the braiding is defined as
 follows:
$$
\Psi_{V,W}:V\o_{t}W\lr
W\o_{t}V ,\ \ \ \ \ v\o_{t}w\longmapsto R^{(2)}\cdot
w\o R^{(1)}\c v.
$$

Let $L$ be a quantum groupoid with antipode $\ul S$. We consider the centralizer of $L_s$ in
$L$, denoted by
$$
C_{L}(L_s)=\{l\in L\mid lx=xl, \, \mbox{for all $x\in L_s$}\},
$$Then we have the following proposition.
\\

 {\bf Proposition 2.1.} Let $(H, \R)$ be a quasitriangular quantum groupoid, $L$ a quantum groupoid and $f$
a quantum groupoid map between $H$ and $L$. Then $C_{L}(L_s)$ is an algebra in $_H\mathcal{M}$, where the
left $H$-module structure on $C_{L}(L_s)$
 is defined by
$$
h\c l=f(h_1)lf(S(h_2)), \ \ \ \mbox{for all} \ \ h\in H, l\in C_{L}(L_s).
$$
The
 multiplication and unit are given by:
 $$
 \m_\R: C_{L}(L_s)\o_t C_{L}(L_s)\lr C_{L}(L_s), \ \ a\o_t b\lm ab,  \ \ \mbox{for all} \ \ a, b\in C_{L}(L_s),
 $$
 $$
 \e_\R: H_t\lr C_{L}(L_s), \ \ \ \  x\lm f(x),  \ \ \mbox{for all} \ \  x\in H_t.
 $$
\\

{\bf Proof.} Firstly, we check that $C_{L}(L_s)$ satisfies the module
conditions. For all $h, g\in H$, and $l\in C_{L}(L_s)$, we have
\begin{eqnarray*}
(gh)\c l=f(g_1)f(h_1)l
f(S(h_2))f(S(g_2))=g\c(h\c l),
\end{eqnarray*}
and $ 1\c
l=f(1_1)lf(S(1_2))=lf(1_1)\ul{S}f(1_2)=l.$

Secondly, we show that $\m_\R$ and $\e_\R$ are morphisms in $_H\mathcal{M}$. In fact,
for all $h\in H$, $a, b\in C_{L}(L_s)$, we have
$$
(h_1\c a)(h_2\c b)=f(h_1)af(\v_s(h_2))bf(S(h_3))=f(h_1)abf(S(h_2))=h\c(ab)
$$
and
$$
\v_t(h)\c1_L=f(h_1S(h_3))f(S(\v_t(h_2)))=f(h_1)f(S(h_2))=h\c1_L.
$$

Thirdly, we prove that $\m_\R$ and $\e_\R$ are well-defined. Obviously, $ab\in C_{L}(L_s)$, $ \mu_\R(a\o_t b)=(1_1\c b)(1_2\c b)=1\c(ab)=ab$ and
\begin{eqnarray*}
f(x)=f(\v_t(x))=f(1_2)\ul{\v}f(1_1x)=\ul{\v}_t(f(x)),
\end{eqnarray*}
Thus we have $f(x)\in L_t\subseteq C_{L}(L_s)$.

Finally, it is easy to see that the
associativity are straightforward. So it remains to verify the unity.
Indeed, for all $z\in H_t, a\in C_{L}(L_s)$, we have
$$
\m_\R\ci(\e_\R\o id)(z\o_t a)=f(z)a=f(1_1z)af(S(1_2))=z\c a,
$$
and
$$
\m_\R\ci(id\o\e_\R)(a\o_t z)=af(z)=f(S^{-1}(1_2))af(1_1z)=a\c z.
$$
This finishes the proof of the proposition.
\\

 {\bf Proposition 2.2.}  Let $(H, \R)$ be a quasitriangular quantum groupoid, $L$ a quantum groupoid and $f$
a quantum groupoid map between $H$ and $L$. Then $C_{L}(L_s)$ is a coalgebra in $_{H}\mathcal{M}$. The
 comultiplication and counit are given by
 $$
\D_\R: C_{L}(L_s)\lr C_{L}(L_s)\o_t C_{L}(L_s),  \ \ \  l\lm l_1f(S(R^{(2)}))\o R^{(1)}\c l_2,
$$
and
$$
\v_\R: C_{L}(L_s)\lr H_t, \ \  l\lm \ul{\v}(f(1_1)l)1_2, \ \ \mbox{for all}
\ \ l\in C_{L}(L_s),
$$
respectively, where $\ul{\v}$ is the counit of $L$.
\\

{\bf Proof:} Firstly, we verify that $\D_\R$ and $\v_\R$ are
well-defined. Obviously, $\v_\R(l)\in H_t$. For
all $l\in C_{L}(L_s)$, we have
\begin{eqnarray*}
\D_\R(l)
&=&l_1f(S(r^{(2)}))f(S(R^{(2)}))\o
f(R^{(1)})l_2f(S(r^{(1)}))\\
&=&f(1_1)l_1f(S(r^{(2)}))f(S(R^{(2)}))\o
f(R^{(1)}1_2)l_2f(S(r^{(1)}))\\
&=&f(1_1)l_1f(S(r^{(2)}))f(S(1_2R^{(2)}))\o
f(R^{(1)})l_2f(S(r^{(1)}))\\
&=&f(1_1)l_1f(S(R^{(2)}))f(S(1_2))\o
R^{(1)}\c l_2\\
&\in& C_{L}(L_s)\o C_{L}(L_s).
\end{eqnarray*}
and
\begin{eqnarray*}
\D_\R(l)
&=&f(1_1)l_1f(S(R^{(2)}r^{(2)}))\o f(R^{(1)})f(1_2)l_2f(S(r^{(1)}))\\
&=&f(1_1)l_1f(S(R^{(2)}r^{(2)}))\o f(R^{(1)})f(S^{-1}(1_3)1_2)l_2f(S(r^{(1)}))\\
&=&f(1_1)l_1f(S(R^{(2)}1_3r^{(2)}))\o f(R^{(1)})f(1_2)l_2f(S(r^{(1)}))\\
&=&f(1_1)l_1f(S(R^{(2)}1_3S(1_4)r^{(2)}))\o f(R^{(1)}1_2)l_2f(S(r^{(1)}))\\
&=&f(1_1)l_1f(S(1_2R^{(2)}S(1_4)r^{(2)}))\o f(1_3R^{(1)})l_2f(S(r^{(1)}))\\
&=&f(1_1)l_1f(S(1_2R^{(2)}r^{(2)}))\o f(1_3)f(R^{(1)})l_2f(S(1_4r^{(1)}))\\
&=&1_1\c (l_1f(S(R^{(2)}))\o 1_2\c (R^{(1)}\c l_2)\\
&\in& C_{L}(L_s)\o_t C_{L}(L_s)
\end{eqnarray*}

Then the coassociativity is same as the setting of Hopf
algebras, then we check the counity as follows:
\begin{eqnarray*}
&&(\v_\R\o id)\ci\D_\R(l)\\
&=&\ul{\v}(f(1_1)l_1f(S(r^{(2)}))f(S(R^{(2)}))f(1_2)f(R^{(1)})l_2f(S(1_3r^{(1)}))\\
&=&\ul{\v}(f(1_1)l_1f(S(S(1_3)r^{(2)}))f(S(R^{(2)}))f(1_2)f(R^{(1)})l_2f(S(r^{(1)}))\\
&=&\ul{\v}(f(1_1)l_1f(S(r^{(2)}))f(S^2(R^{(2)}))f(1_2S(R^{(1)}1_3))l_2f(S(r^{(1)}))\\
&=&\ul{\v}(f(1_1)l_1f(S(r^{(2)}))f(S(R^{(2)}))f(1_2)f(R^{(1)})l_2f(S(r^{(1)}))\\
&=&\ul{\v}_{t}(l_1f(S(r^{(2)}))f(S(R^{(2)})))f(R^{(1)})l_2f(S(r^{(1)}))\\
&=&l_1\ul{\v}_{t}(f(S(R^{(2)})))\ul{S}(l_2)f(R^{(1)}_1)l_3f(S(R^{(1)}_2))\\
&=&l_1\v(1_1S(R^{(2)}))f(1_2)\ul{S}(l_2)f(R^{(1)}_1)l_3f(S(R^{(1)}_2))\\
&=&l_1\v S(R^{(2)}1_2)f(S(1_1))\ul{S}(l_2)(R^{(1)}\c l_3)\\
&=&l_1\v(R^{(2)})f(S(1_1))\ul{S}(l_2)(R^{(1)}S^{-1}(1_2)\c l_3)\\
&=&l_1f(S(1_1))\ul{S}(l_2)(l_3\c 1_2)\\
&=&l_1f(S(1_1))\ul{\v}_s(l_2)f(1_2)=l.
\end{eqnarray*}
Likewise, on the other side we have
\begin{eqnarray*}
&&(id\o\v_\R)\ci\D_\R(l)\\
&=&l_1f(S(r^{(2)}))f(S(R^{(2)}))\ul{\v}_{t}(f(R^{(1)})l_2f(S(r^{(1)})))\\
&=&l_1f(S(r^{(2)}))f(S(R^{(2)}))f(R^{(1)}_1)l_2\ul{\v}_{t}(f(S(r^{(1)})))\ul{S}(l_3)\ul{S}(f(R^{(1)}_2))\\
&=&l_1f(r^{(2)})f(S(x^{(2)}))f(S(R^{(2)}))f(R^{(1)})l_2\ul{\v}(f(1_1r^{(1)}))f(1_2)\ul{S}(l_3)\ul{S}(f(x^{(1)}))\\
&=&l_1f(r^{(2)}1_1)f(S(x^{(2)}))f(S(R^{(2)}))f(R^{(1)})l_2\v(r^{(1)})f(1_2)\ul{S}(l_3)\ul{S}(f(x^{(1)}))\\
&=&l_1f(S(r^{(2)}))f(S(R^{(2)}))f(R^{(1)})\ul{\v}_{t}(l_2)f(S(r^{(1)}))\\
&=&f(1_1)lf(r^{(2)})f(S(R^{(2)}))f(R^{(1)})f(1_2r^{(1)})\\
&=&f(1_1)lf(S(1_2)r^{(2)})f(S(R^{(2)}))f(R^{(1)})f(r^{(1)})\\
&=&lf(S(1_1))f(1_2)=l.
\end{eqnarray*}
Thus we obtain $(\v_\R\o id)\ci\D_\R=id=(id\o\v_\R)\ci\D_\R.$

Finally, we need to prove that $\D_\R$ and $\v_\R$ are morphisms in
$\mathcal{C}$. For all $l\in C_{L}(L_s)$, we do the following calculation:
\begin{eqnarray*}
&&\v_\R(h\c l)=
\ul{\v}(f(1_1)f(h_1)l\ul{S}(f(h_2)))1_2
=\ul{\v}(f(1_1)f(h_1)\ul{\v}_s(f(h_2))l)1_2\\
&=&\ul{\v}(f(1_1h)l)1_2=\ul{\v}(f(h_1)l)\v_t(h_2)=\ul{\v}(f(\v_s(h_1))l)\v_t(h_2)
=\ul{\v}(f(1_1)l)\v_t(h1_2)\\
&=&\ul{\v}(f(1_1)l)h_11_2S(h_2)=h\c\v_\R(l),
\end{eqnarray*}
and $\D_\R(h\c l)=h\c\D_\R(l)$ is clear. This completes the
verification of the coalgebra structure on $C_{L}(L_s)$.
\\

 {\bf Theorem 2.3.} Let $(H, \R)$ be a quasitriangular quantum groupoid, $L$ a quantum groupoid and $f$
a quantum groupoid map between $H$ and $L$. Then $C_{L}(L_s)$ is a bialgebra in $_{H}\mathcal{M}$ with
 structures defined in Proposition 2.1 and Proposition 2.2, and we
 denote $\D_\R(l)= l_{\ul{1}}\o_t l_{\ul{2}}$.
\\

 {\bf Proof.} We now to check that these structure maps obey the
 axioms of a bialgebra in $_{H}\mathcal{M}$. First, we claim that
$(\mu_\R\o \mu_\R)\circ (id\o \Psi_{C_{L}(L_s), C_{L}(L_s)}\o id)\circ(\D_\R\o
\D_\R)=\D_\R\ci\m_\R$, in fact, for all $a, b \in C_{L}(L_s)$,
\begin{eqnarray*}
&&(\mu_\R\o \mu_\R)\circ (id\o \Psi_{C_{L}(L_s), C_{L}(L_s)}\o id)\circ(\D_\R\o
\D_\R)(a\o_t b)\\
&=&a_1f(S(r^{(2)}x^{(2)}))b_{\ul{1}}f(S(R^{(2)}))\o_t(R^{(1)}\c(S(r^{(1)})x^{(1)}\c a_1))b_{\ul{2}}\\
&=&a_1f(S(1_2))b_{\ul{1}}f(S(R^{(2)}))\o_t(R^{(1)}1_1\c a_1)b_{\ul{2}}\\
&=&a_1f(S(1_2))b_{\ul{1}}f(S(R^{(2)}S(1_1)))\o_t(R^{(1)}\c a_1)b_{\ul{2}}\\
&=&a_1f(1_1)b_{\ul{1}}f(S(1_2))f(S(R^{(2)}))\o_t(R^{(1)}\c a_1)b_{\ul{2}}\\
&=&a_1b_1f(S(r^{(2)}))f(S(R^{(2)}))\o_t(R^{(1)}\c a_1)(r^{(1)}\c b_2)\\
&=&a_1b_1f(S(R^{(2)}))\o_t(R^{(1)}_1\c a_1)(R^{(1)}_2\c b_2)\\
&=&a_1b_1f(S(R^{(2)}))\o_tR^{(1)}\c a_2b_2\\
&=&\D_\R\ci\m_\R(a\o_t b).
\end{eqnarray*}
Then we prove $\v_\R \circ \mu_\R= \v_\R \o \v_\R$ as follows:
\begin{eqnarray*}
&&\m_\R(\v_\R(a)\o_t\v_\R(b))
=\ul{\v}(f(1_1)a)\ul{\v}(f(\v_s(1_2)b)1_3\\
&=&\ul{\v}(af(1_1))\ul{\v}(f(1_2)b)1_3
=\ul{\v}(af(1_1)b)1_2=\ul{\v}(f(1_1)ab)1_2\\
&=&\v_\R\ci\m_\R(a\o_t b).
\end{eqnarray*}
Finally, for $1_L\in C_{L}(L_s)$, we have
\begin{eqnarray*}
\D_\R(1_L)
&=&f(1_1)f(r^{(2)})f(S(R^{(2)}))\o_tf(R^{(1)})f(1_2r^{(1)})\\
&=&f(1_1)f(S(1_2)r^{(2)})f(S(R^{(2)}))\o_tf(R^{(1)})f(r^{(1)})\\
&=&f(S(R^{(2)}r^{(2)}))\o_tf(R^{(1)}S(r^{(1)}))\\
&=&f(S(1_1))\o_t f(1_2)
=f(S(1_1))f(1_2)\o_t 1_L\\
&=&1_L\o_t 1_L.
\end{eqnarray*}
This completes our proof.
\\

{\bf Theorem 2.4.} In the situation of Theorem 2.4, $C_{L}(L_s)$ is a Hopf
algebra in $_{H}\mathcal{M}$ with antipode $S_\R$, which is given
by
$$
S_\R: C_{L}(L_s)\lr C_{L}(L_s), \ \  l\lm f(R^{(2)})\ul{S}(R^{(1)}\c l), \ \ l\in C_{L}(L_s).
$$

{\bf Proof.} For convenience we write $S_\R$ in the following equivalent form
$$
S_\R(l)=f(u^{-1})f(S(R^{(2)}))\ul{S}(l)f(R^{(1)}),
$$
where $u^{-1}=R^{(2)}S^2(R^{(1)})$ is the inverse of the element
$u=S(R^{(2)})R^{(1)}$.
First, we assert that $S_\R$ is well-defined. In fact, for all $l\in C_{L}(L_s)$, we have
\begin{eqnarray*}
&&f(1_1)f(u^{-1})f(S(R^{(2)}))S(l)f(R^{(1)})f(S(1_2))\\
&=&f(u^{-1})f(S^2(1_1))f(S(R^{(2)}))S(l)f(R^{(1)})f(S(1_2))\\
&=&f(u^{-1})f(S(R^{(2)}1_2))S(l)f(R^{(1)})f(1_1)\\
&=&f(u^{-1})f(S(R^{(2)}))S(l)f(R^{(1)}S^{-1}(1_2))f(1_1)\\
&=&f(u^{-1})f(S(R^{(2)}))S(l)f(R^{(1)}).
\end{eqnarray*}
Thus, $S_\R(l)\in C_{L}(L_s)$. It is easy to get $S_\R$ is a morphism in
$_{H}\mathcal{M}$. Now we prove the properties of antipode. Then for all
$l\in C_{L}(L_s),$
\begin{eqnarray*}
&&\m_\R\ci(id\o S_\R)\ci\D_\R(l)\\
&=&l_1f(S(R^{(2)}))f(u^{-1})f(S(r^{(2)}))\ul{S}(R^{(1)}\c l_2)f(r^{(1)})\\
&=&l_1f(S(R^{(2)}))f(r^{(2)}u^{-1})\ul{S}(R^{(1)}\c l_2)f(S(r^{(1)}))\\
&=&l_1f(S(r^{(2)}R^{(2)}))f(u^{-1})\ul{S}(R^{(1)}\c l_2)f(S^2(r^{(1)}))\\
&=&l_1f(S(R^{(2)}))f(u^{-1})\ul{S}(f(R^{(1)}_2)l_2f(S(R^{(1)}_3)))f(S^2(R^{(1)}_1))\\
&=&l_1f(S(R^{(2)}))f(u^{-1})f(S^2(R^{(1)}_2))\ul{S}(l_2)f(S(\v_s(R^{(1)}_1)))\\
&=&l_1f(S(R^{(2)}))f(u^{-1})f(S^2(R^{(1)}1_2))\ul{S}(l_2)f(S(1_1))\\
&=&l_1f(R^{(2)})f(u^{-1})f(S(R^{(1)}))\ul{S}(f(S(1_2)))\ul{S}(l_2)\ul{S}(f(1_1))\\
&=&l_1f(u^{-1})f(S^2(R^{(2)}))f(S(R^{(1)}))\ul{S}(l_2)\\
&=&l_1f(u^{-1})f(u)\ul{S}(l_2)=\ul{\v}_t(l)\\
&=&\e_\R\ci\v_\R(l),
\end{eqnarray*}
and
\begin{eqnarray*}
&&\m_\R\ci(S_\R\o id)\ci\D_\R(l)
=f(u^{-1})f(S(R^{(2)}))\ul{S}(l_1)l_2f(R^{(1)})\\
&=&f(u^{-1})f(S(R^{(2)}))\ul{\v}(lf(1_2))f(1_1R^{(1)})
=f(u^{-1})f(S(R^{(2)}1_1))f(R^{(1)})\ul{\v}(lf(1_2))\\
&=&f(u^{-1})f(S^2(1_2))f(u)\ul{\v}(lf(S(1_1)))
=f(1_2)f(u^{-1})f(u)\ul{\v}(lf(S(1_1)))\\
&=&f(1_2)\ul{\v}(l\ul{S}(f(1_1)))
=\ul{\v}(f(1_1)l)f(1_2)=\e_\R\ci\v_\R(l).
\end{eqnarray*}
This completes the proof of the theorem.
\\

{\bf Corollary 2.5.} Let $(H, \R)$ be a quasitriangular quantum groupoid, then $C_H(H_s)$ is a Hopf algebra in the category
$_H\mathcal{M}$, where the left $H$-module structure on $C_H(H_s)$ is given by
$$
h\c g=h_1gS(h_2), \ \ \ \mbox{ for all}\ \  h\in H, g\in C_H(H_s),
$$
and we denote this action of $H$ on $C_H(H_s)$ by $Ad$, i.e., $Ad_{h}(g)=h_1gS(h_2)$. The algebra structure is defined as follows:
$$
\widehat{\m}(g\o_t l)=gl, \ \ \widehat{\e}(x)=x, \ \ \mbox{ for all }\ \ g, l\in C_H(H_s), \ \ x\in H_t.
$$
The comultiplication, counit, and antipode are defined by
$$
\widehat{\D}(g)=g_1S(R^{(2)})\o_t Ad_{R^{(1)}}(g_2), \ \  \widehat{\v}(g)=\v_t(g), \ \  \mbox{ and }\ \
\widehat{S}(h)=R^{(2)}S(R^{(1)}\c h),
$$
respectively.
\\

{\bf Example 2.6.} Let $\mathbb{R}$ be the real field. Set $N=\left\{\left(
                                             \begin{array}{cc}
                                               a & 0 \\
                                               0 & b \\
                                             \end{array}
                                           \right)
\mid a, b\in \mathbb{R} \right\}$, then $N$ is a quantum groupoid (see \cite{WZ10}) with the basis
$$e_1=\left(
        \begin{array}{cc}
          1 & 0 \\
          0 & 0 \\
        \end{array}
      \right), \ \ e_2=\left(
                         \begin{array}{cc}
                           0 & 0 \\
                           0 & 1 \\
                         \end{array}
                       \right).
$$
Its multiplication, unit, comultiplication, counit and antipode are given by
$$
\mu(e_i\o e_j)=\left\{
                \begin{array}{ll}
                  e_i, \ \ i=j \hbox{} \\
                  0, \ \ \ i\neq j \hbox{}
                \end{array}
              \right.; \ \ \
\ \ \ \ \
\eta(1)=:e=e_1+e_2;
$$
$$
\D(e_i)=e_i\o e_i; \ \ \ \ \ \ \v(e_i)=1; \ \ \ \ S(e_i)=e_i; \ \ \   i=1, 2.
$$

Let $\R=e_1\o e_1+e_2\o e_2\in \Delta^{cop}(e)(N\otimes_{k}N)\Delta(e)$.
Then $(N, \R)$ is a quasitriangular quantum groupoid. According to Corollary 2.5, $C_{N}(N_s)=N$ as vector space,
is a Hopf algebra in the category
$_{N}\mathcal{M}$, where the left $N$-module structure is the multiplication on $N$,
and the Hopf algebra structures are defined as follows:
$$
\widehat{\m}=\mu; \ \ \ \widehat{\e}(e_i)=e_i; \ \ \ \widehat{\v}(e_i)=\v_t(e_i)=e_i;
$$
$$
\widehat{\D}(e_i)=e_iS(R^{(2)})\o_t Ad_{R^{(1)}}(e_i)=e_i\o e_i;
$$
$$
\widehat{S}(e_i)=R^{(2)}S(R^{(1)}\c e_i)=e_i; \ \ \ ~i=1, 2. \ \
$$

{\bf Corollary 2.7.} Let $(H, \R)$ be a quasitriangular Hopf algebra, $L$ a Hopf algebra and
$f$ a Hopf algebra map between them. Then $C_L(L_s)=L$ and Theorem 2.4 is the
transmutation theorem in the sense of Majid \cite{M93}.

\section*{3. Quantization approach for quantum groupoids }
\def\theequation{1. \arabic{equation}}
\setcounter{equation} {0} \hskip\parindent

Let $H$ be a quantum gropoid, $F=F^{(1)}\o F^{(2)}\in \D(1)(H\o H)\D^{cop}(1)$ a weak invertible unit 2-cocycle.
In this section, we assume that $H$ is cocommutative and denote the category of left $H$-modules by $_H{\mathcal{M}}_F$,
which module action is $Ad$, braiding and module structure on tensor product are given by
$$
\Phi(m\o_t n)=Ad_{F^{-(1)}F^{(2)}}(n)\o_t Ad_{F^{-(2)}F^{(1)}}(m),
$$
and
$$
Ad_{h}(m\o_t n)=Ad_{F^{-(1)}h_1F^{(1)}}(m)\o_t Ad_{F^{-(2)}h_2F^{(2)}}(n),
$$
respectively, for all $m\in M$, $n\in N$, $M, N\in {_H{\mathcal{M}}_F}$. Then we construct a Hopf algebra in the
category $_H{\mathcal{M}}_F$.
\\

{\bf Proposition 3.1.} Let $H$ be a cocommutative quantum groupoid, $F$ a weak invertible unit 2-cocycle.
Then $C_H(H_s)$ is an algebra in $_H{\mathcal{M}}_F$, which multiplication and unit
are defined by
$$
\c_F: C_H(H_s)\o_t C_H(H_s)\lr C_H(H_s), \ \ \ a\c_F b=Ad_{F^{(1)}}(a)Ad_{F^{(2)}}(b),
$$
and
$$
\e_F: H_t\lr C_H(H_s), \ \ \  \e_F(x)=x,
$$
for all $a, b\in C_H(H_s)$, $x\in H_t$.
\\

{\bf Proof.} First, we show that the multiplication is well-defined. For all $a, b\in C_H(H_s)$, $x\in H_s$, we have
\begin{eqnarray*}
&&\c_F(a\o_t b)
=Ad_{F^{(1)}1_1}(a)Ad_{F^{(2)}1_2}(b)
=Ad_{F^{(1)}}(a)Ad_{F^{(2)}}(b),
\end{eqnarray*}
and by
\begin{eqnarray*}
xAd_{F^{(1)}}(a)&=&xF^{(1)}_1aS(F^{(1)}_2)=x1_1F^{(1)}_1aS(F^{(1)}_2)S(1_2)\\
&=&F^{(1)}_1aS(F^{(1)}_2)x=Ad_{F^{(1)}}(a)x,
\end{eqnarray*}
we have $Ad_{F^{(1)}}(a)\in C_H(H_s)$. Similarly, $Ad_{F^{(2)}}(b)\in C_H(H_s)$, so $a\c_F b\in C_H(H_s)$.

Secondly, we check the associativity and unity. For all $a, b, c\in C_H(H_s)$, we do the following calculation.
\begin{eqnarray*}
(a\c_F b)\c_F c&=&Ad_{F'^{(1)}}(Ad_{F^{(1)}}(a)Ad_{F^{(2)}}(b))Ad_{F'^{(2)}}(c)\\
&=&F'^{(1)}_1Ad_{F^{(1)}}(a)1_1Ad_{F^{(2)}}(b)S(F'^{(1)}_21_2)Ad_{F'^{(2)}}(c)\\
&=&F'^{(1)}_1Ad_{F^{(1)}}(a)\v_s(F'^{(1)}_2)Ad_{F^{(2)}}(b)S(F'^{(1)}_3)Ad_{F'^{(2)}}(c)\\
&=&Ad_{F'^{(1)}_1F^{(1)}}(a)Ad_{F'^{(1)}_2F^{(2)}}(b)Ad_{F'^{(2)}}(c)\\
&=&Ad_{F'^{(1)}}(a)F'^{(2)}_1Ad_{F^{(1)}}(b)S(F'^{(2)}_2)F'^{(2)}_3Ad_{F^{(2)}}(c)S(F'^{(2)}_4)\\
&=&Ad_{F'^{(1)}}(a)F'^{(2)}_1Ad_{F^{(1)}}(b)\v_s(F'^{(2)}_2)Ad_{F^{(2)}}(c)S(F'^{(2)}_3)\\
&=&Ad_{F'^{(1)}}(a)F'^{(2)}_1\v_s(F'^{(2)}_2)Ad_{F^{(1)}}(b)Ad_{F^{(2)}}(c)S(F'^{(2)}_3)\\
&=&Ad_{F'^{(1)}}(a)F'^{(2)}_1Ad_{F^{(1)}}(b)Ad_{F^{(2)}}(c)S(F'^{(2)}_2)\\
&=&a\c_F(b\c_F c),
\end{eqnarray*}
and for all $z\in H_t$,  we have
\begin{eqnarray*}
&&\c_F(id\o \e_F)(a\o_t z)=Ad_{F^{(1)}}(a)Ad_{F^{(2)}}(z)\\
&=&Ad_{F^{(1)}}(a)\v_t(F^{(2)}z)=Ad_{zF^{(1)}}(a)\v_t(F^{(2)})=Ad_{z1_1}(a)1_2\\
&=&z_1a\v_s(z_2)=\v_t(z)a=\v(1_1z)a1_2=az=a(Ad_z(1)),
\end{eqnarray*}
and
\begin{eqnarray*}
&&\c_F(\e_F\o id)(z\o_t a)=Ad_{F^{(1)}}(z)Ad_{F^{(2)}}(a)=\v_t(F^{(1)}z)Ad_{F^{(2)}}(a)\\
&=&\v_t(F^{(1)})Ad_{F^{(2)}S^{-1}(z)}(a)=S(1_1)Ad_{1_2S^{-1}(z)}(a)=S(1_1)1_2azS(1_3)\\
&=&a\v_t(z)=1_2a\v(1_1z)=za=Ad_z(a).
\end{eqnarray*}

Finally, it is clear that $\e_F$ is a morphism in $_H\mathcal{M}_F$, so we only need to prove that $\c_F$ is a morphism in $_H\mathcal{M}_F$. In fact, we have
\begin{eqnarray*}
&&\c_F(Ad_h(a\o_t b))=
Ad_{F'^{(1)}F^{-(1)}h_1F^{(1)}}(a)Ad_{F'^{(2)}F^{-(2)}h_2F^{(2)}}(b)\\
&=&Ad_{1_1h_1F^{(1)}}(a)Ad_{1_2h_2F^{(2)}}(b)
=h_1Ad_{F(1)}(a)\v_s(h_2)Ad_{F(2)}(b)S(h_3)\\
&=&h_1\v_s(h_2)Ad_{F(1)}(a)Ad_{F(2)}(b)S(h_3)
=Ad_{h}(a\o_F b).
\end{eqnarray*}
This completes the proof.
\\

{\bf Proposition 3.2.} Let $H$ be a cocommutative quantum groupoid, $F$ a weak invertible unit 2-cocycle.
Then $C_H(H_s)$ is a coalgebra in $_H{\mathcal{M}}_F$, which comultiplication and counit
are given by
$$
\D_F: C_H(H_s)\lr C_H(H_s)\o_t C_H(H_s), \ \ \ \D_F(a)=Ad_{F^{-(1)}}(a_1)\o_t Ad_{F^{-(2)}}(a_2),
$$
and
$$
\v_F: C_H(H_s)\lr H_t, \ \ \ \ \ \v_F(a)=\v_t(a),
$$
for all $a\in C_H(H_s)$.
\\

{\bf Proof. } First, we prove that $\D_F$ and $\v_F$ are well-defined. Obviously, $\v_F$ is well-defined.
Now we check as follows,
for all $a\in C_H(H_s)$, we have
\begin{eqnarray*}
Ad_{F^{-(1)}}(a_1)\o_t Ad_{F^{-(2)}}(a_2)=Ad_{1_1F^{-(1)}}(a_1)\o Ad_{1_2F^{-(2)}}(a_2)=Ad_{F^{-(1)}}(a_1)\o Ad_{F^{-(2)}}(a_2),
\end{eqnarray*}
Here the second equality uses the cocommutativity.

And, for all $x\in H_s$, we have
\begin{eqnarray*}
xAd_{F^{-(1)}}(a_1)=x1_1F^{-(1)}_1a_1S(F^{-(1)}_2)S(1_2)=Ad_{F^{-(1)}}(a_1)x,
\end{eqnarray*}
i.e., $Ad_{F^{-(1)}}(a_1)\in C_H(H_s)$. Similarly, $Ad_{F^{-(2)}}(a_2)\in C_H(H_s)$.
Then $Ad_{F^{-(1)}}(a_1)\o_t Ad_{F^{-(2)}}(a_2)\in C_H(H_s)\o_t C_H(H_s)$.

Secondly, the coassociativity is straightforward, then we verify the counity. In fact, we have
\begin{eqnarray*}
&&(\v_F\o id)\D_F(a)
=\v_t(F^{-(1)}_1a_1S(F^{-(1)}_2))\o_t Ad_{F^{-(2)}}(a_2)\\
&=&\v_t(F^{-(1)}_1a_1\v_s(F^{-(1)}_2))\o_t Ad_{F^{-(2)}}(a_2)=\v_t(F^{-(1)}1_1a_1S(1_2))\o_t Ad_{F^{-(2)}}(a_2)\\
&=&\v_t(F^{-(1)}a_1S(1_2))\o_t Ad_{F^{-(2)}S(1_1)}(a_2)=\v_t(F^{-(1)}a_11_1)\o_t Ad_{F^{-(2)}1_2}(a_2)\\
&=&\v_t(F^{-(1)}a_1S(1_1))\o_t Ad_{F^{-(2)}}(1_2a_2S(1_3))
=\v_t(F^{-(1)}S(1'_1))\o_t Ad_{F^{-(2)}}(1_11'_2aS(1_2))\\
&=&\v_t(F^{-(1)}S(1_1))\o_t Ad_{F^{-(2)}1_2}(a)
=Ad_{\v_t(F^{-(1)})F^{-(2)}}(a)=a
\end{eqnarray*}
and
\begin{eqnarray*}
&&(id\o\v_F)\D_F(a)=Ad_{F^{-(1)}}(a_1)\o_t\v_t (F^{-(2)}_1a_2\v_s(F^{-(2)}_2))\\
&=&Ad_{F^{-(1)}}(a_1)\o_t\v_t (F^{-(2)}1_1a_2S(1_2))=Ad_{F^{-(1)}}(a_1)\o_t\v_t (F^{-(2)}1_2a_2S(1_1))\\
&=&Ad_{F^{-(1)}S^{-1}(1_2)}(a_1)\o_t\v_t (F^{-(2)}a_21_1)=Ad_{F^{-(1)}}(1_1a_1S(1_2))\o_t\v_t (F^{-(2)}a_2S(1_3))\\
&=&Ad_{F^{-(1)}}(1_11'_1aS(1_2))\o_t\v_t (F^{-(2)}1'_2)=Ad_{F^{-(1)}}(1_2a)\o_t\v_t (F^{-(2)}1_1)\\
&=&Ad_{F^{-(1)}1_2}(a)\o_t\v_t (F^{-(2)}1_1)
=Ad_{S^{-1}(\v_t(F^{-(2)}))F^{-(1)}}(a)=a.
\end{eqnarray*}

Finally, we prove that $\D_F$ and $\v_F$ are morphisms in $_H\mathcal{M}_F$. We compute
\begin{eqnarray*}
\D_F(Ad_h(a))&=&Ad_{F^{-(1)}}(h_1a_1S(h_4))\o_t Ad_{F^{-(2)}}(h_2a_2S(h_3))\\
&=&Ad_{F^{-(1)}h_1}(a_1)\o_t Ad_{F^{-(2)}h_2}(a_2)\\
&=&Ad_{F^{-(1)}h_1F^{(1)}F'^{-(1)}}(a_1)\o_t Ad_{F^{-(2)}h_2F^{(2)}F'^{-(2)}}(a_2)\\
&=&Ad_h(\D_F(a))
\end{eqnarray*}
and
\begin{eqnarray*}
\v_F(Ad_h(a))&=&\v_t(h_1aS(h_2))=\v_t(h_1\v_s(h_2)a)=\v_t(ha)\\
&=&h_1\v_t(a)S(h_2)=Ad_h(\v_F(a))
\end{eqnarray*}
as required. This completes the proof of the proposition.
\\

{\bf Theorem 3.3.} Let $H$ be a cocommutative quantum groupoid, $F$ a weak invertible unit 2-cocycle.
Then $C_H(H_s)$ is a Hopf algebra in $_H{\mathcal{M}}_F$, which antipode
$S_F:  C_H(H_s)\lr C_H(H_s)$ is given by $S_F(a)=S(a)$, for all $a\in C_H(H_s)$.
\\

{\bf Proof. } In order to show that $C_H(H_s)$ is a Hopf algebra in $_H{\mathcal{M}}_F$,
 we prove that $\D_F$ and $\v_F$ are algebra maps in $_H{\mathcal{M}}_F$ first.
 It is easy to get the following equality.
 \begin{eqnarray*}
 &&Y^{(1)}F^{-(1)}\o Y^{(2)}X^{-(1)}F^{(2)}Y^{-(1)}\o X^{(1)}X^{-(2)}F^{(1)}F^{-(2)}\o X^{(2)}Y^{-(2)}\\
 &~&\qquad\qquad=F^{-(1)}_1F^{(1)}_1\o F^{-(1)}_2F^{(2)}_1\o F^{-(2)}_1F^{(1)}_2\o F^{-(2)}_2F^{(2)}_2,
 \end{eqnarray*}
 where $F=Y=X$, and $F^{-(1)}=Y^{-(1)}=X^{-(1)}$.
 For all $a, b\in C_H(H_s)$, by the equality, we have
 \begin{eqnarray*}
 &&\D_F(a\o_F b)=\D_F(Ad_{F^{(1)}}(a)Ad_{F^{(2)}}(b))\\
&=&Ad_{F^{-(1)}}(F^{(1)}_1a_1S(F^{(1)}_4)F^{(2)}_1b_1S(F^{(2)}_4))
 \o_t Ad_{F^{-(2)}}(F^{(1)}_2a_2S(F^{(1)}_3)F^{(2)}_2b_2S(F^{(2)}_3))\\
&=&Ad_{F^{-(1)}}(Ad_{F^{(1)}_1}(a_1)Ad_{F^{(2)}_1}(b_1))
 \o_t Ad_{F^{-(2)}}(Ad_{F^{(1)}_2}(a_2)Ad_{F^{(2)}_2}(b_2))\\
&=&F^{-(1)}_1Ad_{F^{(1)}_1}(a_1)Ad_{F^{(2)}_1}(b_1)S(F^{-(1)}_2)
 \o_t Ad_{F^{-(2)}}(Ad_{F^{(1)}_2}(a_2)Ad_{F^{(2)}_2}(b_2))\\
&=&F^{-(1)}_1\v_s(F^{-(1)}_2)Ad_{F^{(1)}_1}(a_1)Ad_{F^{(2)}_1}(b_1)S(F^{-(1)}_3)
 \o_t Ad_{F^{-(2)}}(Ad_{F^{(1)}_2}(a_2)Ad_{F^{(2)}_2}(b_2))\\
&=&F^{-(1)}_1Ad_{F^{(1)}_1}(a_1)\v_s(F^{-(1)}_2)Ad_{F^{(2)}_1}(b_1)S(F^{-(1)}_3)
 \o_t Ad_{F^{-(2)}}(Ad_{F^{(1)}_2}(a_2)Ad_{F^{(2)}_2}(b_2))\\
&=&Ad_{F^{-(1)}_1F^{(1)}_1}(a_1)Ad_{F^{-(1)}_2F^{(2)}_1}(b_1)
 \o_t Ad_{F^{-(2)}_1F^{(1)}_2}(a_2)Ad_{F^{-(2)}_2F^{(2)}_2}(b_2)\\
&=&Ad_{Y^{(1)}F^{-(1)}}(a_1)Ad_{Y^{(2)}X^{-(1)}F^{(2)}Y^{-(1)}}(b_1)
\o_t Ad_{X^{(1)}X^{-(2)}F^{(1)}F^{-(2)}}(a_2)Ad_{X^{(2)}Y^{-(2)}}(b_2)\\
&=&Ad_{F^{-(1)}}(a_1)\c_F Ad_{X^{-(1)}F^{(2)}Y^{-(1)}}(b_1)
\o_t Ad_{X^{-(2)}F^{(1)}F^{-(2)}}(a_2)\c_F Ad_{Y^{-(2)}}(b_2)\\
&=&(Ad_{F^{-(1)}}(a_1)\o_t Ad_{F^{-(2)}}(a_2))(Ad_{Y^{-(1)}}(b_1)\o_t Ad_{Y^{-(2)}}(b_2))\\
&=&\D_F(a)\D_F(b),
 \end{eqnarray*}
and
\begin{eqnarray*}
&&\D_F(1)=Ad_{F^{-(1)}}(1_1)\o_t Ad_{F^{-(2)}}(1_2)=Ad_{F^{-(1)}}(1_1)\o_t Ad_{F^{-(2)}1_2}(1)\\
&=&Ad_{F^{-(1)}}(1_2)\o_t Ad_{F^{-(2)}1_1}(1)=Ad_{F^{-(1)}1_2}(1)\o_t Ad_{F^{-(2)}1_1}(1)\\
&=&\v_t(F^{-(1)})\o_t \v_t(F^{-(2)})=1\o_t \v(1_1F^{-(1)})1_2\v_t(F^{-(2)})\\
&=&1\o_t \v(F^{-(1)})\v_t(F^{-(2)})=1\o_t 1.
\end{eqnarray*}

Thus, we obtain $\D_F$ is an algebra map in $_H{\mathcal{M}}_F$. Now, we show that
$\v_F$ is an algebra map as follows:
\begin{eqnarray*}
\v_F(a\c_F b)&=&\v_t(Ad_{F^{(1)}}(a)Ad_{F^{(2)}}(b))=\v_t(Ad_{F^{(1)}}(a)F^{(2)}_1b\v_s(F^{(2)}_2))\\
&=&\v_t(Ad_{F^{(1)}}(a)F^{(2)}b)=\v(1_1Ad_{F^{(1)}}(a)F^{(2)}b)1_2\\
&=&\v(\v_s(F^{(1)}_1)aS(F^{(1)}_2)1_1F^{(2)}b)1_2=\v(aS(F^{(1)})1_1F^{(2)}b)1_2\\
&=&\v(aS(F^{(1)})S(1_1))\v_t(1_2F^{(2)}b)=\v(aS(F^{(1)}))\v_t(F^{(2)}b)\\
&=&\v(a\v_s(F^{(1)}))\v_t(F^{(2)}b)=\v(aS(1_2))\v_t(1_1b)\\
&=&\v(1_2a)\v_t(1_1b)=\v_t(\v_t(a)b)\\
&=&\v_F(a)\v_F(b).
\end{eqnarray*}

So we obtain $C_H(H_s)$ is a bialgebra in $_H{\mathcal{M}}_F$. Then, we want to
verify the properties of antipode. Obviously, $S_F$ is well-defined, and $S_F(Ad_h(a))=Ad_h(S_F(a))$ holds for
all $h\in H$. Next, we do the following calculations:
\begin{eqnarray*}
\c_F(id\o_t S_F)\D_F(a)
&=&Ad_{F^{(1)}F^{-(1)}}(a_1)Ad_{F^{(2)}}(S^2(F^{-(2)}_2)S(a_2)S(F^{-(2)}_1))\\
&=&Ad_{F^{(1)}F^{-(1)}}(a_1)Ad_{F^{(2)}F^{-(2)}}(S(a_2))=Ad_{1_1}(a_1)Ad_{1_2}(S(a_2))\\
&=&1_1a_1S(1_3)S(1_2a_2)=1_2a_1S(1_1)S(1_3a_2)=\v_t(1_2a1_1)\\
&=&\v_t(1_2aS(1_1))=\v_t(a)=\v_F(a)
\end{eqnarray*}
and
\begin{eqnarray*}
&&\c_F(S_F\o id)\D_F(a)
=Ad_{F^{(1)}}(S(Ad_{F^{-(1)}}(a_1))Ad_{F^{(2)}}(Ad_{F^{-(2)}}(a_2))\\
&=&Ad_{F^{(1)}F^{-(1)}}(S(a_1))Ad_{F^{(2)}F^{-(2)}}(a_2)
=1_1S(a_1)\v_s(1_2)S^2(a_2)S(1_3)\\
&=&1_1S(a_2)\v_s(1_2)S^2(a_1)S(1_3)=\v_t(S(a))=\v_t(a)=\v_F(a).
\end{eqnarray*}
This completes our proof.
\\

In this paper, we denote the Hopf algebra we obtained in Theorem 3.3 by $C_H(H_s)_F$.
\\

{\bf Example 3.4.} Let $N$ be the quantum groupoid in Example 2.6, $F=F^{(1)}\o F^{(2)}=e_1\o e_1+e_2\o e_2\in \Delta(e)(N\otimes_{k}N)\Delta^{cop}(e)$
 a weak invertible unit 2-cocycle. According to Theorem 3.3, $C_{N}(N_s)_F=N$ as vector space,
is a Hopf algebra in the category
$_{N}\mathcal{M}_F$, where the Hopf algebra structures on $C_{N}(N_s)_F$ are defined as follows:

$$
e_i\c_F e_j=Ad_{F^{(1)}}(e_i)Ad_{F^{(2)}}(e_j)=\left\{
                                                 \begin{array}{ll}
                                                   e_i, \ \ i=j \hbox{} \\
                                                   0, \ \ \ i\neq j \hbox{}
                                                 \end{array}
                                               \right.
; \ \ \ \e_F(e_i)=e_i;
$$
$$
\D_F(e_i)=Ad_{F^{-(1)}}(e_i)\o_t Ad_{F^{-(2)}}(e_i)=e_i\o e_i; \ \ \ \v_F(e_i)=\v_t(e_i)=e_i;
$$
$$
S_F(e_i)=S(e_i)=e_i; \ \ \ i=1, 2.
$$

{\bf Corollary 3.5.} Let $H$ be a cocommutative Hopf algebra, $F$ a cocycle.
Then $C_H(H_s)=H$, and $H_F$ is an $S$-Hopf algebra in the sense of \cite{GRZ92}.
The multiplication and comultiplication are given by
$$
h\c_F g=Ad_{F^{(1)}}(h)Ad_{F^{(2)}}(g)\ \ \mbox{ and }\ \
\D_F(h)=Ad_{F^{-(1)}}(h_1)\o Ad_{F^{-(2)}}(h_2),
$$
for all $h, g\in H_F$, respectively. The unit, counit and antipode of $H_F$ coincide with those of
$H$.

\section*{4. Braided groups obtained by weak twisting thoery }
\def\theequation{1. \arabic{equation}}
\setcounter{equation} {0} \hskip\parindent

Let $(H, \R)$ be a quasitriangular quantum groupoid, $F$ a weak invertible unit 2-cocycle, then by the
result of \cite{C08}
there is a new quasitriangular quantum groupoid $(\widetilde{H}, \widetilde{\R})$, defined by the same
multiplication, unit and counit, and for all $h\in\wt H$,
$$
\widetilde\D(h)=F^{-1}\D(h)F, \ \ \ \widetilde\R=F^{-1}_{21}\R F, \ \ \ \widetilde S(h)=vS(h)v^{-1},
$$
where $v=F^{(-1)}S(F^{-(2)})$ and $v^{-1}=S(F^{(1)})F^{(2)}$.
\\

Applying Corollary 2.5, we obtain that $C_{\widetilde H}(\widetilde H_s)$ is a Hopf algebra in
the category $_{\widetilde H}\mathcal{M}$ of the left $\wt H$-modules, which module action is given by
$$
Ad_{\wt h}(g)=h_{\wt 1}g\wt S(h_{\wt 2}), \ \ \mbox{ for all } \ \
h\in \wt H, \ \ g\in C_{\widetilde H}(\widetilde H_s),
$$
where $\wt\D(h)=h_{\wt 1}\o h_{\wt 2}$, $h\in\widetilde H$.
\\

{\bf Lemma 4.1.} The element $v^{-1}=S(F^{(1)})F^{(2)}$ obeys the following equality:
$$
\D(v^{-1})=((S\o S)(F^{-1}_{21}))(v^{-1}\o v^{-1})F^{-1}.
$$

{\bf Proof. } We calculate as follows:
\begin{eqnarray*}
\D(v^{-1})
&=&S(F^{(1)}_2)F^{(2)}_1\o S(F^{(1)}_1)F^{(2)}_2=S(F^{(1)}_2)F^{(2)}_11_1\o S(F^{(1)}_1)F^{(2)}_21_2\\
&=&S(F^{(1)}_2)F^{(2)}_1X^{(1)}X^{-(1)}\o S(F^{(1)}_1)F^{(2)}_2X^{(2)}X^{-(2)}\\
&=&S(X^{(1)}_2)\v_s(F^{(1)}_2)X^{(2)}X^{-(1)}\o S(X^{(1)}_1)S(F^{(1)}_1)F^{(2)}X^{-(2)}\\
&=&S(1_2X^{(1)}_2)X^{(2)}X^{-(1)}\o S(1_1X^{(1)}_1)S(F^{(1)})F^{(2)}X^{-(2)}\\
&=&S(X^{(1)}_21_2)X^{(2)}X^{-(1)}\o S(X^{(1)}_11_1)v^{-1}X^{-(2)}\\
&=&S(X^{(1)}_2F^{(2)}F^{-(2)})X^{(2)}X^{-(1)}\o S(X^{(1)}_1F^{(1)}F^{(-1)})v^{-1}X^{-(2)}\\
&=&S(F^{-(2)})S(X^{(2)}_1F^{(1)})X^{(2)}_2F^{(2)}X^{-(1)}\o S(F^{(-1)})S(X^{(1)})v^{-1}X^{-(2)}\\
&=&S(F^{-(2)})S(F^{(1)})\v_s(X^{(2)})F^{(2)}X^{-(1)}\o S(F^{(-1)})S(X^{(1)})v^{-1}X^{-(2)}\\
&=&S(\v_s(X^{(2)})F^{-(2)})S(F^{(1)})F^{(2)}X^{-(1)}\o S(F^{(-1)})S(X^{(1)})v^{-1}X^{-(2)}\\
&=&S(F^{-(2)})S(F^{(1)})F^{(2)}X^{-(1)}\o S(S^{-1}(\v_s(X^{(2)}))F^{(-1)})S(X^{(1)})v^{-1}X^{-(2)}\\
&=&S(F^{-(2)})v^{-1}X^{-(1)}\o S(F^{(-1)})\v_s(X^{(2)})S(X^{(1)})v^{-1}X^{-(2)}\\
&=&S(F^{-(2)})v^{-1}X^{-(1)}\o S(F^{(-1)})v^{-1}X^{-(2)}.
\end{eqnarray*}
This finishes our proof.
\\

{\bf Proposition 4.2.} Let $H$ be cocommutative and $F$ a weak invertible unit 2-cocycle. Then the category
$_{\widetilde H}\mathcal{M}$ can be identified with the category $_H{\mathcal{M}}_F$.
\\

{\bf Proof.} Following \cite{GM94}, it is straightforward.
\\

According to Section 3, we obtain that $C_H(H_s)_F$ is a Hopf algebra in the category ${_H{\mathcal{M}}_F}$,
then we have the following lemma.
\\

{\bf Lemma 4.3.} There is an isomorphism
 $\a: C_H(H_s)_F\lr C_{\widetilde H}(\widetilde H_s)$ of
$\widetilde H$-modules given by $\a(a)=Ad_{F^{(1)}}(a)F^{(2)}$, for all $a\in C_H(H_s)_F$,
and $C_H(H_s)_F\in {_H{\mathcal{M}}_F}$, $C_{\wt H}(\wt H_s)\in {_{\wt H}\mathcal{M}}$,
where the structures of $C_H(H_s)_F$ and $C_{\wt H}(\wt H_s)$ are defined as before.
\\

{\bf Proof. } First of all, we compute another useful form of $\a$ as follows.
\begin{eqnarray*}
&&\a(a)=F^{(1)}_1aS(F^{(1)}_2)F^{(2)}=F^{(1)}_1X^{(1)}X^{-(1)}aS(X^{(2)}X^{-(2)})S(F^{(1)}_2)F^{(2)}\\
&=&F^{(1)}X^{-(1)}aS(X^{-(2)})S(F^{(2)}_1X^{(1)})F^{(2)}_2X^{(2)}
=F^{(1)}X^{-(1)}aS(X^{-(2)})S(X^{(1)})\v_s(F^{(2)})X^{(2)}\\
&=&1_2X^{-(1)}aS(X^{-(2)})S(X^{(1)})1_1X^{(2)}=1_2X^{-(1)}aS(X^{-(2)})S(X^{(1)}1_1)X^{(2)}\\
&=&F^{-(1)}aS(F^{-(2)})v^{-1}.
\end{eqnarray*}
Then we easy to see that $\a(a)\in C_{\widetilde H}(\widetilde H_s)$.
In fact, for all $a\in C_H(H_s)_F$, we have
\begin{eqnarray*}
\widetilde {Ad}_1(\a(a))&=&1_2\a(a)vS(1_1)v^{-1}=1_2F^{-(1)}aS(F^{-(2)})v^{-1}vS(1_1)v^{-1}\\
&=&F^{-(1)}aS(F^{-(2)})v^{-1}=\a(a).
\end{eqnarray*}

And the equivalent form of $\a$ imply that $\a$ is invertible with $\a^{-1}(a)=F^{(1)}avS(F^{(2)})$. Obviously,
$\a^{-1}(a)\in C_H(H_s)_F$, and we have
\begin{eqnarray*}
\a\circ\a^{-1}(a)&=&\a(F^{(1)}avS(F^{(2)}))=F^{-(1)}F^{(1)}avS(F^{(2)})S(F^{-(2)})v^{-1}\\
&=&1_2avS(1_1)v^{-1}=1_2aF^{-(1)}S(F^{-(2)})S(1_1)v^{-1}\\
&=&1_1aS^{-1}(1_2)=a.
\end{eqnarray*}
$\a^{-1}\circ\a(a)=a$ is straightforward. The rest proof of the lemma is the same as one in the Hopf algebra case.
This finishes the proof.
\\

{\bf Lemma 4.4.} Let $H$ be cocommutative and $F$ a weak invertible unit 2-cocycle.
The map $\a: C_H(H_s)_F\lr C_{\widetilde H}(\widetilde H_s)$ in Lemma 4.3 is
an algebra map. \\

{\bf Proof.} For all $a, b\in C_H(H_s)_F$, we calculate
\begin{eqnarray*}
\a(a\c_F b)&=&Ad_{F^{(1)}}(a)F^{(2)}_1Ad_{X^{(1)}}(b)S(F^{(2)}_2)F^{(2)}_3X^{(2)}\\
&=&Ad_{F^{(1)}}(a)F^{(2)}_1Ad_{X^{(1)}}(b)\v_s(F^{(2)}_2)X^{(2)}\\
&=&Ad_{F^{(1)}}(a)F^{(2)}1_1Ad_{X^{(1)}}(b)S(1_2)X^{(2)}\\
&=&Ad_{F^{(1)}}(a)F^{(2)}Ad_{X^{(1)}}(b)X^{(2)}\\
&=&\a(a)\a(b),
\end{eqnarray*}
as required.
\\

For the remainder of the section we assume that the quansitriangular structure of $H$ is $\R=\D^{cop}(1)\D(1)$.
\\

{\bf Lemma 4.5.} Let $H$ be cocommutative and $F$ a weak invertible unit 2-cocycle.
 Then the map $\a: C_H(H_s)_F\lr C_{\widetilde H}(\widetilde H_s)$ in Lemma 4.3 is
a coalgebra map.
\\

{\bf Proof.} For all $a \in C_H(H_s)_F$, firstly, we compute the expression of $\D(\a(a))$.
\begin{eqnarray*}
&&\D(\a(a))=F^{-(1)}_1a_1S(F^{-(1)}_2)(v^{-1})_1\o F^{-(1)}_2a_2S(F^{-(2)}_1)(v^{-1})_2\\
&=&F^{-(1)}_1a_1S(F^{-(1)}_2)S(F'^{-(2)})v^{-1}X^{-(1)}\o F^{-(1)}_2a_2S(F^{-(2)}_1)S(F'^{-(1)})v^{-1}X^{-(2)}\\
&=&(F^{-(1)}_1F'^{-(1)}_1a_1S(F'^{-(2)})\o F^{-(1)}_2F'^{-(1)}_2a_2S(F^{-(2)}F'^{-(1)}_3))(v^{-1}\o v^{-1})F^{-1}.
\end{eqnarray*}
Then, we have
\begin{eqnarray*}
\widetilde {\D}(\a(a))=F^{-(1)}\a(a)_1F^{(1)}\o F^{-(2)}\a(a)_2F^{(2)}.
\end{eqnarray*}

Next, applying Corollary 2.6 to $(\wt H, \wt\R)$, we obtain
\begin{eqnarray*}
\widehat{\wt\D}(\a(a))&=&F^{-(1)}\a(a)_1F^{(1)}\widetilde R^{(2)}\wt S(\wt r^{(2)})
\o_t \wt r^{(1)}F^{-(2)}\a(a)_2F^{(2)}\wt R^{(1)}\\
&=&F^{-(1)}\a(a)_1F^{(1)}Y^{-(1)}Y^{(2)}vS(X^{-(1)}X^{(2)})v^{-1}
\o_t X^{-(2)}X^{(1)}F^{-(2)}\\
&&\a(a)_2F^{(2)}Y^{-(2)}Y^{(1)}\\
&=&F^{-(1)}\a(a)_11_1Y^{(2)}vS(F^{(2)})S(Y^{-(1)})v^{-1}
\o_t Y^{-(2)}F^{(1)}F^{-(2)}\a(a)_21_2Y^{(1)}\\
&=&F^{-(1)}\a(a)_2Y^{(2)}vS(F^{(2)})S(Y^{-(1)})v^{-1}
\o_t Y^{-(2)}F^{(1)}F^{-(2)}\a(a)_1Y^{(1)}\\
&=&F^{-(1)}X^{-(1)}_2F'^{-(1)}_2a_2S(F'^{-(1)}_3)S(X^{-(2)})v^{-1}Y'^{(2)}Y^{(2)}vS(F^{(2)})S(Y^{-(1)})v^{-1}\\
&&\o_t Y^{-(2)}F^{(1)}F^{-(2)}X^{-(1)}_1F'^{-(1)}_1a_1S(F'^{-(2)})v^{-1}Y'^{(1)}Y^{(1)}\\
&=&F^{-(1)}X^{-(1)}_2Ad_{F'^{-(1)}_2}(a_2)S(X^{-(2)})v^{-1}1_1X'^{-(1)}S(X'^{-(2)})S(F^{(2)})S(Y^{-(1)})\\
&&v^{-1}\o_t Y^{-(2)}F^{(1)}F^{-(2)}X^{-(1)}_1F'^{-(1)}_1a_1S(F'^{-(2)})S(Y^{(1)})Y^{(2)}1_2\\
&=&F^{-(1)}X^{-(1)}_2Ad_{F'^{-(1)}_2}(a_2)S(X^{-(2)})v^{-1}X'^{-(1)}S(X'^{-(2)}1_1)S(F^{(2)})S(Y^{-(1)})\\
&&v^{-1}\o_t Y^{-(2)}F^{(1)}F^{-(2)}X^{-(1)}_1F'^{-(1)}_1a_1S(F'^{-(2)})S(1_2Y^{(1)})Y^{(2)}\\
&=&F^{-(1)}X^{-(1)}_2Ad_{F'^{-(1)}_2}(a_2)S(X^{-(2)})v^{-1}X'^{-(1)}S(1_1)S(X'^{-(2)})S(F^{(2)})\\
&&S(Y^{-(1)})v^{-1}\o_t Y^{-(2)}F^{(1)}F^{-(2)}X^{-(1)}_1F'^{-(1)}_1a_1S(F'^{-(2)})S(1_2)S(Y^{(1)})Y^{(2)}\\
&=&F^{-(1)}X^{-(1)}_2Ad_{F'^{-(1)}_2}(a_2)S(X^{-(2)})v^{-1}X'^{-(1)}S(S(1_1)X'^{-(2)})S(F^{(2)})\\
&&S(Y^{-(1)})v^{-1}\o_t Y^{-(2)}F^{(1)}F^{-(2)}X^{-(1)}_1F'^{-(1)}_1a_1S(F'^{-(2)}S(1_2))S(Y^{(1)})Y^{(2)}\\
&=&F^{-(1)}X^{-(1)}_2Ad_{F'^{-(1)}_2}(a_2)S(X^{-(2)})S(1_2)S(F^{(2)})S(Y^{-(1)})v^{-1}\o_t Y^{-(2)}F^{(1)}\\
&&F^{-(2)}X^{-(1)}_1F'^{-(1)}_1a_1S(F'^{-(2)})S(1_1)v^{-1}\\
&=&F^{-(1)}X^{-(1)}_1Ad_{F'^{-(1)}_1}(a_1)S(F^{(2)}1_2X^{-(2)})S(Y^{-(1)})v^{-1}\o_t Y^{-(2)}F^{(1)}F^{-(2)}\\
&&X^{-(1)}_2F'^{-(1)}_2a_2S(1_1F'^{-(2)})v^{-1}\\
&=&X^{-(1)}Ad_{F'^{-(1)}_1}(a_1)S(X^{-(2)}_2)S(F^{-(2)})S(F^{(2)}1_2)S(Y^{-(1)})v^{-1}\o_t Y^{-(2)}F^{(1)}\\
&&F^{-(1)}X^{-(2)}_1F'^{-(1)}_2a_2S(1_1F'^{-(2)})v^{-1}\\
&=&X^{-(1)}Ad_{F'^{-(1)}_1}(a_1)S(X^{-(2)}_1)S(F^{-(2)})S(F^{(2)})S(Y^{-(1)})v^{-1}\o_t Y^{-(2)}1_2F^{(1)}\\
&&F^{-(1)}X^{-(2)}_2F'^{-(1)}_2a_2S(1_1F'^{-(2)})v^{-1}\\
&=&F^{-(1)}Ad_{F'^{-(1)}_1}(a_1)S(F^{-(2)}_1)S(1_2)S(Y^{-(1)})v^{-1}\o_t Y^{-(2)}1'_21_1F^{-(2)}_2F'^{-(1)}_2\\
&&a_2S(1'_1F'^{-(2)})v^{-1}\\
&=&F^{-(1)}1_2Ad_{F'^{-(1)}_1}(a_1)S(Y^{-(1)}F^{-(2)}_1)v^{-1}
\o_t Y^{-(2)}F^{-(2)}_2F'^{-(1)}_2a_2S(1_1F'^{-(2)})v^{-1}\\
&=&Y^{-(1)}F^{-(1)}_11_2Ad_{F'^{-(1)}_1}(a_1)S(Y^{-(2)}F^{-(1)}_2)v^{-1}
\o_t F^{-(2)}F'^{-(1)}_2a_2S(1_1F'^{-(2)})v^{-1}\\
&=&Y^{-(1)}F^{-(1)}_11_2Ad_{F'^{-(1)}_1}(a_1)S(1_3)S(F^{-(1)}_2)S(Y^{-(2)})v^{-1}\o_t F^{-(2)}F'^{-(1)}_2\\
&&a_2S(1_1F'^{-(2)})v^{-1}\\
&=&Y^{-(1)}Ad_{F^{-(1)}1_2F'^{-(1)}_1}(a_1)S(Y^{-(2)})v^{-1}\o_t F^{-(2)}F'^{-(1)}_2a_2S(1_1F'^{-(2)})v^{-1}\\
&=&Y^{-(1)}Ad_{F^{-(1)}F'^{-(1)}_1}(a_1)S(Y^{-(2)})v^{-1}\o_t 1_2F^{-(2)}F'^{-(1)}_2a_2S(1_1F'^{-(2)})v^{-1}\\
&=&Y^{-(1)}Ad_{F'^{-(1)}}(a_1)S(Y^{-(2)})v^{-1}\o_t 1_2F^{-(1)}Ad_{F'^{-(2)}}(a_2)S(1_1F^{(2)})v^{-1}\\
&=&Y^{-(1)}Ad_{F'^{-(1)}}(a_1)S(Y^{-(2)})v^{-1}\o_t F^{-(1)}Ad_{F'^{-(2)}}(a_2)S(F^{(2)})v^{-1}\\
&=&(\a\o\a)\D_F(a).
\end{eqnarray*}

Furthermore, we have
\begin{eqnarray*}
\widehat{\wt{\v}}\circ\a(a)
&=&\v_t(F^{-(1)}aS(F^{-(2)})S(F^{(1)})F^{(2)})\\
&=&\v_t(F^{-(1)}F^{(1)}_1aS(F^{(1)}_2)S(F^{-(2)}_1)F^{-(2)}_2F^{(2)})\\
&=&\v_t(F^{-(1)}Ad_{F^{(1)}}(a)\v_s(F^{-(2)})F^{(2)})\\
&=&\v_t(1_1Ad_{F^{(1)}}(a)S(1_2)F^{(2)})\\
&=&\v_t(Ad_{F^{(1)}}(a)\v_t(F^{(2)}))\\
&=&\v_t(1_1a\v_s(1_2))=\v_F(a).
\end{eqnarray*}
This completes the proof.
\\

{\bf Lemma 4.6.} Let $H$ be cocommutative and $F$ a weak invertible unit 2-cocycle.
Then the antipode $S_F$ on $C_H(H_s)_F$ and the antipode $\widehat{\wt{S}}$ on $C_{\wt H}({\wt H}_s)$
satisfy the following condition:
$$
\widehat{\wt{S}}\circ \a=\a\circ S_F.
$$

{\bf Proof.} In order to verification, we calculate as follows:
\begin{eqnarray*}
&&\a^{-1}\circ\widehat{\wt{S}}\circ \a(a)\\
&=&F^{(1)}Y^{-(1)}Y^{(2)}X^{-(1)}S(X^{-(2)})
S(X^{(2)})Ad_{X^{(1)}Y^{-(2)}Y^{(1)}}(S(a))S(F^{(2)})\\
&=&F^{(1)}Y^{-(1)}X^{-(1)}_2Y^{(2)}_1S(X^{-(2)}Y^{(2)}_2)
S(X^{(2)})Ad_{X^{(1)}Y^{-(2)}X^{-(1)}_1Y^{(1)}}(S(a))S(F^{(2)})\\
&=&F^{(1)}Y^{-(1)}X^{-(1)}_2\v_t(Y^{(2)})S(X^{-(2)})
S(X^{(2)})Ad_{X^{(1)}Y^{-(2)}X^{-(1)}_1Y^{(1)}}(S(a))S(F^{(2)})\\
&=&F^{(1)}Y^{-(1)}X^{-(1)}_21_2S(X^{-(2)})
S(X^{(2)})Ad_{X^{(1)}Y^{-(2)}X^{-(1)}_11_1}(S(a))S(F^{(2)})\\
&=&F^{(1)}Y^{-(1)}X^{-(1)}_1S(X^{(2)}X^{-(2)})Ad_{X^{(1)}Y^{-(2)}X^{-(1)}_2}(S(a))S(F^{(2)})\\
&=&F^{(1)}X^{-(1)}S(X^{(2)}Y^{-(2)}X^{-(2)}_2)Ad_{X^{(1)}Y^{-(1)}X^{-(2)}_1}(S(a))S(F^{(2)})\\
&=&F^{(1)}X^{-(1)}S(1_2X^{-(2)}_2)Ad_{1_1X^{-(2)}_1}(S(a))S(F^{(2)})\\
&=&F^{(1)}Y^{-(1)}S(Y^{-(2)}_2)Ad_{Y^{-(2)}_1}(S(a))S(F^{(2)})\\
&=&F^{(1)}Y^{-(1)}S(Y^{-(2)}_3)Y^{-(2)}_1S(a)S(Y^{-(2)}_2)S(F^{(2)})\\
&=&F^{(1)}Y^{-(1)}S(Y^{-(2)}_1)Y^{-(2)}_2S(a)S(Y^{-(2)}_3)S(F^{(2)})\\
&=&F^{(1)}Y^{-(1)}\v_s(Y^{-(2)}_1)S(a)S(Y^{-(2)}_2)S(F^{(2)})\\
&=&F^{(1)}Y^{-(1)}1_1S(a)S(1_2)S(Y^{-(2)})S(F^{(2)})\\
&=&1_1S(a)S(1_2)=S(a)=S_F(a).
\end{eqnarray*}
Thus we have $\widehat{\wt{S}}\circ \a=\a\circ S_F$. This finishes the proof.
\\

The following isomorphism theorem is the main result of this section.
\\

{\bf Theorem 4.7.} Let $H$ be cocommutative and $F$ a weak invertible unit 2-cocycle. Assume that
$C_H(H_s)_F$ be the Hopf algebra with the structures as shown in Theorem 3.3 viewed
as an object in the category $_{\wt H}\mathcal{M}$ of $\wt H$-modules. In this category,
there is an isomorphism of Hopf algebras
$$
\a: C_H(H_s)_F\cong C_{\widetilde H}(\widetilde H_s).
$$

{\bf Example 4.8.} Let $N$ be the quantum groupoid in Example 2.6, $\R=F=e_1\o e_1+e_2\o e_2$ the quasitriangular
structure and weak invertible unit 2-cocycle. Then we construct a new quasitriangular quantum groupoid $(\widetilde N, \widetilde \R)$
with the same multiplication, unit and counit, and
$$
\widetilde\D(e_i)=F^{-1}\D(e_i)F=e_i\o e_i, \ \ \ \widetilde\R=F^{-1}_{21}\R F=e_1\o e_1+e_2\o e_2,
$$
$$
\widetilde S(e_i)=F^{(-1)}S(F^{-(2)})S(e_i)S(F^{(1)})F^{(2)}=e_i, \ \ \ \ i=1, 2.
$$
So $(N, \R)=(\widetilde N, \widetilde \R)$ as quasitrianlar quantum groupoids. By Example 2.6 and Example 3.4, we obtain that
$$
C_N(N_s)_F=C_{\widetilde N}(\widetilde N_s),
$$
as objects in the category $_{\widetilde N}\mathcal{M}$.
\\

{\bf Corollary 4.9.} Let $H$ be a cocommutative Hopf algebra, and $F$ a cocycle. Then the Theorem 4.7 is
Theorem 2.8 in \cite{GM94}.
\\

\begin{center}
 {\bf ACKNOWLEDGEMENT}
\end{center}

 The work was partially supported by the NNSF of China (No.11326063), and Qing Lan Project.

\end{document}